\documentclass[11pt]{article}
\usepackage[english]{babel}
\usepackage{amsfonts}
\usepackage[dvips]{graphicx}

\begin{document}

\title{\bf Exact solution of a differential problem in analytical fluid dynamics: use of Airy's functions}

\date{June 2006}

\author{\bf Gianluca Argentini \\
\normalsize gianluca.argentini@gmail.com \\
\normalsize gianluca.argentini@riellogroup.com \\
\textit{Research \& Development Department}\\
\textit{Riello Burners}, 37048 San Pietro di Legnago (Verona), Italy}

\maketitle

\begin{abstract}
Treating a boundary value problem in analytical fluid dynamics, translation of 2D steady Navier-Stokes equations to ordinary differential form leads to a second order equation of Riccati type. In the case of a compressible fluid with constant kinematic viscosity along streamlines, it is possible to find an exact solution of the differential problem by rational combination of Airy's functions and their derivatives.
\end{abstract}

\newtheorem{theo}{Theorem}
\newtheorem{prop}{Proposition}
\newtheorem{corol}{Corollary}
\newtheorem{deft}{Definition}

\section{Ordinary differential form of Navier-Stokes equations}

Suppose to have a 2D steady flow of a fluid. Let

\begin{equation}\label{parameterization}
	\Phi : s \longmapsto (\phi_1(s),\phi_2(s))=(x,y)
\end{equation}

\noindent an admissible parameterization (\cite{diffg}) for each streamline, with $\Phi: [a,b] \rightarrow {\mathbb R}^2$ for some suitable values $a$ and $b$ such that $\Phi(a)$ is the initial point of the streamline, $\Phi(b)$ is the end point in the considered geometrical domain. Then, if $\mathbf{v}=(v_1,v_2)$ is the flow velocity field, by definition of streamline there is a scalar function $f=f(s)$ such that

\begin{equation}\label{velocityParam}
	(v_1,v_2) = f(s)(\dot{\phi}_1(s),\dot{\phi}_2(s))
\end{equation}

\noindent where $\dot{\phi}$ is the derivative of $\phi$.\\

\noindent \textit{Remark}. The author is developing a more general model of parameterization for Navier-Stokes equations, in the case 3D too, under less restrictive hypothesis than those formulated in this work.\\

\noindent Suppose that $f(s) \neq 0$, $\dot{\phi}_1(s) \neq 0$ for all $s$ along the streamline, that is $v_1 \neq 0$, and let $\mathbf{u}=\mathbf{v} \circ \phi$. Then

\begin{prop}
Along a streamline $s \longmapsto \phi(s)$ the following relation holds: 
\begin{equation}\label{nonLinearTerm}
	\mathbf{v} \cdot \nabla v_1 = \frac{1}{\dot{\phi}_1}u_1\dot{u}_1
\end{equation}
\end{prop}

\noindent \textit{Dim}. Using chain rule we have

\begin{equation}
	\dot{u}_1 = \frac{\partial v_1}{\partial x}\dot{\phi}_1 + \frac{\partial v_1}{\partial y}\dot{\phi}_2
\end{equation}

\noindent But from (\ref{velocityParam})

\begin{equation}
	\dot{\phi}_2 = \frac{v_2}{v_1}\dot{\phi}_1
\end{equation}

\noindent therefore

\begin{equation}
	\dot{u}_1 = \frac{\dot{\phi}_1}{v_1}\left[v_1\frac{\partial v_1}{\partial x} + v_2\frac{\partial v_1}{\partial y}\right]
\end{equation}

\noindent Then (\ref{nonLinearTerm}) follows from the fact that along the streamline $u_1(s)=v_1(\phi(s))$. $\square$\\

\noindent Now let $\phi$ invertible, that is $\exists$ $g: \phi([a,b]) \rightarrow [a,b]$ such that $s=g(\phi(x,y))$ when $(x,y)=\phi(s)$. Then

\begin{prop}
Along a streamline the following relation holds: 
\begin{equation}
	\Delta v_1 = (\nabla g \cdot \nabla g) \ddot{u}_1 + (\Delta g) \dot{u}_1
\end{equation}
\end{prop}

\noindent \textit{Dim}. From definition, we have $\mathbf{v} = \mathbf{u} \circ g$, therefore $\nabla v_1 = \dot{u}_1 \nabla g$. Then

\begin{equation}
	\frac{\partial^2 v_1}{\partial x^2} = \ddot{u}_1 \left( \frac{\partial g}{\partial x} \right)^2 + \dot{u}_1 \left( \frac{\partial^2 g}{\partial x^2} \right)
\end{equation}

\noindent and 

\begin{equation}
	\frac{\partial^2 v_1}{\partial y^2} = \ddot{u}_1 \left( \frac{\partial g}{\partial y} \right)^2 + \dot{u}_1 \left( \frac{\partial^2 g}{\partial y^2} \right)
\end{equation}

\noindent therefore the thesis follows by summing the two previous expressions. $\square$\\

\noindent The following Proposition holds for a scalar function $p=p(x,y)$:

\begin{prop}
Let $q = p \circ \phi$. Along a streamline: 
\begin{equation}
	\Delta p = \dot{q}(\Delta g)
\end{equation}
\end{prop}

\noindent \textit{Dim}. From $q = p \circ \phi$ follows $p = q \circ g$. Then, using chain rule,

\begin{equation}
	\frac{\partial p}{\partial x} = \frac{dq}{ds}\frac{\partial g}{\partial x}
\end{equation}

\noindent and analogous equation for $y$-differentiation. $\square$\\

\noindent Now let $\rho$ the density, $\mu$ the dynamic viscosity and $\mathbf{f}$ the body force per unit volume of fluid. Then the Navier-Stokes equations for the steady flow are (see \cite{batchelor})

\begin{equation}\label{NS}
	\rho (\mathbf{v} \cdot \nabla \mathbf{v}) = -\nabla p + \rho \mathbf{f} + \mu \Delta \mathbf{v}
\end{equation}

\noindent Identifying not derived functions $h(x,y)$ with their composition $h \circ \phi$ and using simple substitutions, from previous Propositions and (\ref{velocityParam}) we can state the \textit{ordinary differential form} of Navier-Stokes equations:

\begin{eqnarray}\label{ordinaryNS}
	\rho \frac{1}{\dot{\phi}_1}u_1\dot{u}_1 & = & \rho f_1 - \dot{q} \frac{\partial g}{\partial x} +\mu (\nabla g \cdot \nabla g) \ddot{u}_1 + \mu(\Delta g) \dot{u}_1 \\
	u_2 & = & \frac{\dot{\phi}_2}{\dot{\phi}_1}u_1 \nonumber
\end{eqnarray}

\noindent Note that, if the streamlines $\phi=(\phi_1,\phi_2)$ are known, from equation (\ref{ordinaryNS}) it could be possible calculate $u_1$. But what about the physical meaning of this solution? We make now some considerations about conservation of mass.\\

\section{About continuity equation}

In the case of incompressible flow, the differential form of this equation is simply $\nabla \cdot \mathbf{v} = 0$  (\cite{batchelor} or \cite{madani}). In this case, from $\mathbf{v} = \mathbf{u} \circ g$, follows that $\partial_x v_1 = \dot{u}_1 \partial_x g$ and $\partial_y v_2 = \dot{u}_2 \partial_y g$, therefore along a streamline the continuity equation is

\begin{equation}\label{continuityTemp}
	\dot{\mathbf{u}} \cdot \nabla g = 0
\end{equation}
\noindent But it is possible, in the hypothesis $\dot{\phi}_1 \neq 0$ along the streamline, to rewrite this equation using only $u_1$:

\begin{prop}
If $\dot{\phi}_1(s) \neq 0$ $\forall \hspace{0.1cm} s$, the continuity equation for incompressible flow is 
\begin{equation}\label{continuity}
	\dot{u}_1 + \left[ \ddot{\phi}_2 - \frac{\dot{\phi}_2}{\dot{\phi}_1}\ddot{\phi}_1 \right]\frac{\partial g}{\partial y}u_1 = 0
\end{equation}
\end{prop}
	
\noindent \textit{Dim}. Differentiating the identity $g(\phi(s))=s$ respect to variable $s$, we have $\nabla g \cdot \dot{\phi}=1$. Then

\begin{equation}
	\partial_x g = \frac{1}{\dot{\phi}_1} - \frac{u_2}{u_1} \partial_y g
\end{equation}

\noindent Substituting this formula in the continuity equation (\ref{continuityTemp}) and using $u_1 \dot{\phi}_2 = u_2 \dot{\phi}_1$, we obtain the new identity

\begin{equation}\label{eqTemp1}
	\dot{u}_1 + [ \dot{\phi}_1\dot{u}_2 - \dot{\phi}_2\dot{u}_1 ]\partial_y g = 0
\end{equation}

\noindent Differentiating the relation $u_2 = \frac{\dot{\phi}_2u_1}{\dot{\phi}_1}$, follows that

\begin{equation}
	\dot{u}_2 = \frac{\ddot{\phi}_2u_1+\dot{\phi}_2\dot{u}_1}{\dot{\phi}_1} - \frac{\ddot{\phi}_1\dot{\phi}_2u_1}{\dot{\phi}_1^2}
\end{equation}

\noindent from which the relation (\ref{continuity}) is obtained eliminating $\dot{u}_2$ in (\ref{eqTemp1}). $\square$\\

\noindent \textit{Remark}. Suppose the streamlines are straight lines expressed by the parameterization $\phi_1(s)=s$, $\phi_2=k$, with $k$ real constants. Then $\dot{\phi}_1=1$ and $\dot{\phi}_2=0$, and from (\ref{continuity}) the continuity equation is simply $\dot{u}_1=0$. Then the velocity field is constant along a streamline, as known for incompressible rectilinear flow (see e.g. \cite{madani}).\\

\noindent As we would investigate existence of non trivial solutions to Navier-Stokes equations (\ref{ordinaryNS}) in the simple case $x=\phi_1(s)=s$, we have $g_x=1$ and $g_y=0$, therefore in the incompressible case from previous Proposition follows that $u_1$ is constant. For a more interesting and realistic solution we assume that flow is steady but with spatially variable density.

\section{A steady compressible flow}

Consider a 2D steady flow where streamlines are parameterizable by the following expressions

\begin{equation}
	x = \phi_1(s) = s, \hspace{0.5cm} y = \phi_2(s) = \phi_2(x)
\end{equation}

\noindent where $\phi_2$ is invertible on an interval $[0,L]$, that is for each streamline exists a function $g: \phi([0,L]) \rightarrow [a,L]$ such that $s=g(x,y)=x$. The function $\phi_2$ and $g$ can depend on the single streamline. Note that $\dot{\phi}_1=1$, $\nabla g = (1,0)$ and $\Delta g = 0$, therefore from (\ref{ordinaryNS}) the Navier-Stokes equations along a streamline become

\begin{equation}\label{problemNS1}
	\rho u_1 \dot{u}_1 = \rho f_1 - \dot{q} + \mu \ddot{u}_1, \hspace{0.5cm} u_2 = \dot{\phi}_2 u_1
\end{equation}

\noindent Suppose that $x$-component $f_1$ of body force $\mathbf{f}$ is constant. Along a streamline we can made the realistic hypothesis that $\rho(s) \neq 0$ $\forall s$; dividing the two members by $\rho$ we obtain

\begin{equation}\label{problemNS2}
	u_1 \dot{u}_1 = f_1 - \frac{\dot{q}}{\rho} + \nu \ddot{u}_1
\end{equation}

\noindent where $\nu=\frac{\mu}{\rho}$ is the kinematic viscosity (see \cite{batchelor}). At this point we suppose that, along a single streamline, $\nu$ and the quantity $\frac{\dot{q}}{\rho}$ are constant, with values depending on streamline. The latter, equivalent to $p=k_1\int\rho+k_2$, can be view as a constitutive equation about the fluid.

\section{An analytical resolution}

In this section we try to find a general exact solution of the non linear differential equation (\ref{problemNS2})$_1$. Integrating on the independent variable $s$ and then dividing by $\nu$, we have

\begin{equation}\label{resolution1}
	\dot{u}_1 = \frac{1}{2\nu}u_1^2 + \frac{1}{\nu}\left(\frac{\dot{q}}{\rho}-f_1\right)s + \frac{c}{\nu}
\end{equation}

\noindent where $c$ is an arbitrary constant. This is a first order differential equation of \textit{Riccati} type (see \cite{boyce} or \cite{goldstein}). Applying the transformation

\begin{equation}\label{transformation}
	u_1 = -2\nu\frac{\dot{z}}{z}
\end{equation}

\noindent we translate previous equation into (see \cite{boyce})

\begin{equation}
	\ddot{z} + \frac{1}{2\nu^2} \left[ \left(\frac{\dot{q}}{\rho}-f_1\right)s + c \right] z = 0
\end{equation}

\noindent This is a form of \textit{Airy}'s type equation (see \cite{goldstein} and \cite{diffEqMath}) and its general solution is (see \cite{diffEqMath})

\begin{equation}
	z(s) = c_1 Ai(t) + c_2 Bi(t)
\end{equation}

\noindent where 

\begin{eqnarray}
	t & = & -\frac{as+b}{(-a)^{\frac{2}{3}}} \\
	a & = & \frac{1}{2\nu^2}\left(\frac{\dot{q}}{\rho}-f_1\right) \\
	b & = & \frac{c}{2\nu^2}	
\end{eqnarray}

\noindent while $Ai(t)$ and $Bi(t)$ are \textit{Airy's functions} (see e.g. \cite{temme}), linearly independent solutions of \textit{Airy's equation} $\ddot{y}-ty=0$ which appears in optics and quantum mechanics phenomena.

\noindent Using transformation (\ref{transformation}), the exact solution of Navier-Stokes equation (\ref{resolution1}) along a streamline is

\begin{equation}\label{solutionNS}
	u_1 = 2\nu (-a)^{\frac{1}{3}} \frac{c_1 Ai'(t) + c_2 Bi'(t)}{c_1 Ai(t) + c_2 Bi(t)}
\end{equation}

\noindent where $Ai'(t)$ and $Bi'(t)$ are the derivative of $Ai$ and $Bi$. \\
Note that this solution has physical meaning only if $a<0$. Assuming that inflow zone is at $s=0$, usually pressure drops down towards outflow, so that we can assume $\dot{q}<0$. Also, we assume $f_1>0$, hence the condition $a<0$ is satisfied.\\
Also, note that with our assumptions $u_1$ doesn't depend on the streamline parametric representation $(x,y)=\phi(s)$; the flow velocity field depends on $\phi$ through $u_2$ component by relation $u_2=\dot{\phi}_2u_1$.

\section{A boundary value problem}

The general solution (\ref{solutionNS}) depends on three constants of integration: $c_1$, $c_2$ and, from (\ref{resolution1}), $c$. Suppose we want to resolve a boundary value problem for equation (\ref{problemNS2}) with $u_1(0)=u_{10}, \dot{u}_1(0)=\dot{u}_{10}$ and $u_1(L)=u_{1L}$. The first two conditions, using (\ref{resolution1}) at $s=0$, give the value of $c$. Then, noting that for $t \in \mathbb{R}$ Airy's functions $Ai(t)$ and $Bi(t)$ have real values (see \cite{temme}), the other two constants $c_1$ and $c_2$ can be found solving the algebraic system $u_1(0)=u_{10}$, $u_1(L)=u_{1L}$.


\begin{thebibliography}{9}

\bibitem{batchelor} Batchelor, G.K. {\it An introduction to fluid dynamics}, Cambridge Mathematical Library, Cambridge University Press (2000)

\bibitem{boyce} Boyce, W.E. and DiPrima, R.C. {\it Elementary Differential Equations}, 6th edition, Wiley, (1996)

\bibitem{goldstein} Braun, W.H. and Goldstein, M.E. {\it Advanced Methods for the solution of Differential Equations}, NASA SP-316, Washington (1973)

\bibitem{diffEqMath} Coombes, K.; Hunt, B.; Lipsman, R.; Osborn, J. and Stuck, G. {\it Differential Equations with Mathematica}, 2nd edition, John Wiley \& Sons, Inc. (1998)

\bibitem{diffg} Lipschutz, M. {\it Differential Geometry}, MacGraw-Hill, New York (1969)

\bibitem{madani} Malek-Madani, R. {\it Advanced Engineering Mathematics}, Addison-Wesley, (1998)

\bibitem{temme} Temme, N.M. {\it Special functions: an introduction to the classical functions of mathematical physics}, J.Wiley, New York (1996)

\end{thebibliography}
\end{document}